# Zeev Schuss' legacy: an unexpected path in applied mathematics

D. Holcman

Ecole Normale Superieure

**Abstract:**

Z. Schuss (1937-2018) was an applied mathematician, with several contributions in asymptotic, stochastic processes, PDEs, modeling and signal processing. He is well known for his original approach to the activation escape problem, based on WKB and boundary layer analysis, summarized in six books published by Springer. The text summarizes his intellectual approach in science, views and personal path across the XX century, WWII and his personal contribution to academia.

**Introduction**

July 29, 2018, passed away the very great Zeev Schuss, Emeritus Professor of Applied Mathematics at Tel Aviv University. There are certain personalities that stand above; above national and international institutions, lobbies, prizes and political forces of all kinds. Not because these people are divine, but because are simply human, hard-working and intellectually superior, not as tested by primitive IQ tests, but by creativity, inspiration and a projecting vision.

Zeev was born in Poland 1937, graduated in composition, conducting, and theory from the Academy of Music in Tel Aviv 1963, graduated in mathematics 1965 and got his PhD in mathematics from Northwestern 1970. He became professor at TAU and served chairman of applied mathematics 1993-1995. He published over 200 papers in pure and applied math, chemistry, physics, engineering, and biology. He wrote 6 books on applied math, published by Springer and Wiley. He supervised tens of MSc and PhD students in the above disciplines; quite a few of them hold positions in prestigious institutions across the world (Israel, USA, Europe, and China).

**Early career**

It took time for Zeev to define for himself what applied mathematics means [15]. After his PhD, he realized that he had to escape from the alienation of proving existence, uniqueness of PDEs and all exercises of finding lower and upper estimates. After his PhD with A. Friedman, Zeev had enough confidence in himself, to move orthogonally to the ambient orthodoxy. The starting point was a class that he took from Henry McKean about stochastic processes and the asymptotics of differential equations that he developed in discussions with B. Matkowsky. Then he discovered the review (published in 1943) of Chandrasekar in the Review of Modern Physics about thermal activation escape from an attractor. He was able to obtain a formal asymptotic formula as a solution of the exit problem in n dimensions. After he gave a talk at the Courant Institute in the 70s, the news spread at the speed of light across the USA [15]. Although not considered as a true mathematical result by the radicals and purists, Zeev understood that this was a new result that nobody could get before and carried in that direction full speed for the next 40 years. It was the mathematics that he wanted to develop and apply to the sciences, engineering, technology and more. His purpose was to analyze, compute, simulate and discover new features, like experimentalists make discovery by looking under the microscope, except that for him, mathematics played the role of the microscope. The goal of applied mathematics for Zeev was to discover, not formal proofs, but to find new physical mechanisms through modeling, new explanations from computations, and to define the new computations needed and to carry them out as accurately as possible.

One of the most robust and efficient tools for finding physical laws are closed formulas, obtained by the method of asymptotic approximation to solutions of PDEs, although the initial formulation is often stochastic processes, probability, statistical physics, chemistry or mechanics. These formulas are particularly valuable for resolving singular behavior, where practically infinite computer time is needed to explore a modest fraction of the parameter space. Such formulas deal precisely with manipulating infinities and are thus very relevant in understanding the refined properties of the studied systems [1-6].

Designing fast and efficient simulations is also a key to success [3]. In this context, finding a proof is not the most urgent mission and is certainly not the goal. Asymptotic approaches have been used successfully more than once, by Poincare for studying divergent series or by Ramanujan in summation of series in number theory or by Chandrasekar in many statistical physics problems. Zeev made several significant contributions applying asymptotics to rare events, such as thermally activated escape from an attractor in physics, chemistry, and in loss of lock in signal tracking [1-3]. Another innovation is the formula for the narrow escape time in molecular and cellular biophysics [4-6].

A large portion of Zeev's intellectual activity was dedicated to applied mathematics of science and engineering, to design stochastic simulations, to analyze stochastic processes, and to develop mathematical biophysics and theoretical biology. In that context, Zeev developed with all of us (his collaborators, friends and students) new tools to analyze data about selectivity of ionic channels (how ions are selected in a channel pore, the size of which can be just a few atomic diameters) and to develop the narrow escape theory (escape of a stochastic particle from a narrow window), that inspired many scientific communities of physics, biophysics, computational biology but also the fourth episode of the TV Fargo Series [16].

**Fighting ideologies**

The reformulation of asymptotics and its formulas into propositions and theorems was easy and he agreed to do this exercise for his last four books [2-6]. But, living as an applied mathematician did not come for free: Zeev always complained about some fanatic pure mathematicians, usually proud of their ignorance of all sciences (biology, chemistry, physics, etc.) and in any mathematical applications that prevented, by pure ideology [10], the development of a leading department of applied math at TAU, which to this day still suffers, although efforts are currently made to change direction. By cutting themselves from Science and technology, the influence of traditional mathematical departments has drastically declined, reduced nowadays to the beauty context of the Fields medal, in which nobody really knows the reasons of somebody deserving such a prize. This is equivalent to an orthodox religion because it is made on belief, no proofs.

Any department that will make the effort to reconnect or connect to data, to medicine, to biology, to internet, etc., will certainly advance. It is never too late. Nowadays individual pioneers are scattered around the globe. Zeev designed new education programs [11] at Tel Aviv University (TAU), to train students in multiple fields of math-physics-and applied computer science.

Zeev like Gauss, Riemann, Wiener, Kolmogorov, liked all mathematics and their applications, but he was very picky in choosing a good problem. Not all problems are equal, not all are interesting, not all are relevant. The choice is key.

**A second career:**

Zeev did not take retirement seriously and kept working full speed for the past 15 years. There is here a deep lesson to learn: when you are good with your students and colleagues, in due time, they will be good to you and you will be entitled to a second career with them, to stay connected to the new science, changing drastically every 7 years. This is what he did, in addition to writing 5 books, published by Springer. A typical day would start at 11pm until 2pm, before going to lunch, where Zeev met his friends, retired professors from all sciences, to discuss any subject including science.

**Music until his melanoma stopped him.**

Zeev started his education by music and kept playing piano until he lost the ability to use his right arm following multiple surgeries of a recurrent melanoma. Zeev fought his severe melanoma by designing a recipe for himself that allowed him to be a singular patient: he made pictures of the skin a few times a month, and compared the color changes of the skin, while developing a sense of touch to detect irregularities in the softness of the tissues. This improvised method allowed him to survive 20 years, while the average time is usually one year for a severe melanoma. No article was written about his case, because you cannot write an article with n=1 in medicine.

Zeev played piano in church in Chicago, while getting other secondary jobs, as taxi driver to provide the complement of salary he needed to support his family when he was a PhD student. For a meeting about experimentalist meeting theoreticians that we organized in 2007 [7], he agreed to play a recital dedicated to Chopin [8].

**His multiple stories with the Teheran Children**

In 2017 in Akko we organized a symposium, where we illustrated the unexpected path of Zeev, paved with obstacles, but also opportunities of the developing sciences by using the geographical journey of the Teheran children [12-15], where Zeev participated both as the youngest child, but also as a defender in a trial that lasted 10 years. To establish the historical truth about ultimately what we call in modern history, the price of being a victim, Zeev decided to participate as a witness, but not as a direct participant in the trial. He said that the demands for repair were nothing: can few tens of thousands shekel compensate for an education that the children never received that impacted their entire life? The money that was supposed to serve this cause, was hijacked and used for something else: to build roads, public services, etc…. During such a too long trial, many of the survivors died, but the truth came up.

Briefly [12], the "Tehran Children" is the story of Jewish children orphan that escaped from Poland to Russia after Germany conquered it in September 1939 with a journey of 4 years, finally arriving in Palestine in 1943, with much funding from a Christian American woman who lived in Iowa*. In 1939, following the German invasion of Poland, 300,000 refugees, including Jews, poles, children, etc., were moving east. In 1940, they were sent to gulags in Siberia in cattle cars. The deportees lived under very difficult conditions, where many of the children died or became orphans in this period. Following the attack on the U.S.S.R by Germany, a new era started for the refugees, where they were released from gulags, and emigrated towards the Asian territories (Uzbekistan, Tajikistan). In 1941, after General Wladyslaw Andres was released from a Russian prison, the Soviet authorities agreed to allow 24,000 Polish citizens to move with the Andres army, including around 1000 Jewish children. 25,000 Polish soldiers were sent to Iran, to strengthen the British armies in the Middle East. Thirty-three thousand soldiers left, 11,000 citizens, 3,000 children, of which about 1,000 were orphaned Jewish children. In January 1943, the children moved to Afhaz and then to the Iranian port of Bender Shapur, where they embarked on the S.S. Dunera, headed to Karachi. From Karachi they embarked on another ship, the Neurolia, which sailed to Suez, Egypt. Then they crossed the Sinai Desert by train, and finally after an odyssey of four years of agony, they arrived in Atlit in Israel ("Palestine" at that time), on February 18th, 1943.

The Jewish agency distributed the children upon their arrival and Zeev Schuss was sent to a kibbutz. Finally, his mother found him and brought him back to Poland, where he could get an education, in contrast with the other kids that stayed in Israel. After the war was over, under the pressure of the USA, West Germany had to pay Israel for the costs of "resettling so great a number of uprooted and destitute Jewish refugees" and to compensate individuals [13].

However, the Teheran children did not receive any compensation and the money that was supposed to help them was used for other purposes. After retirement, many of the Teheran children understood that something went wrong in the process of compensation and started to investigate, until they realized that they were once again victims. In 2006-2007, a group of retired Teheran children initiated a prosecution against the state of Israel, but not much was achieved until Schuss discovered the strategic mistake: after a historian was called upon by the court to give his testimony about the historical context, Zeev found the turning argument: he discovered that "individual compensation" had been mis-translated into "collective compensation". Obviously, this was a deliberated political maneuver to hijack money. This argument affected the defense severely to the point that the judge understood that the money was not intended to build roads but to resettle the refugees, because you do not resettle people on roads, but in houses. More examples were given, and the judge decided to give immediate compensation to the Teheran children. The Appeal came, and the State won and then it went to the supreme court, which confirmed the Appeal. This decision was merely intended to defend the interest of the state: Nobody can attack the State and certainly not victims, especially during the Ben Gurion time. But the State did not ask that the money given to the Teheran children should be returned, proving in a way that the children were right. Once again, there is no glory to be a victim. Zeev fought all the time to get rid of the victim thinking, as he explained in a talk given at Yad Vashem.

**His students and the school**

I met Schuss in 2000 and since then we have been working nonstop, writing tens of papers, 2 books, reviews, opinions, Youtube videos about subjects we thought were important. We usually worked until 8pm every day. Zeev ended the day by saying that there is always another day.

After moving to Israel as a postdoc, I realized that it was a good time to get connected to real science, biology and biophysics, which was impossible to do in France' traditional sciences at that time. I looked at scholarly articles about singular perturbations, and found most of them boring, written for machine not for human, until I found the papers of Schuss where motivations, explanations were given, with deep computations. I understood immediately that this guy knew what he was talking about and was interested in real life problems and their solutions. To my good fortune, my colleague at the Weizmann Institute, M. Margaliot—now Professor at Tel Aviv University, introduced me to Zeev. We met the first day for 3 hours and the second for another 3 hours. This was the beginning of our close collaboration.

In general, Zeev worked hours with every student or colleague to define the question. What was that the correct question to look at, why should we work on that question, why this question is important and for what: does it solve a specific problem in science? It took us usually months before we converged to the correct question, its mathematical formulation and about the same amount of time to resolve it and to design and/or run simulations.

Based on this principle of analysis, this attitude toward science, technology, medicine, etc., is Zeev's legacy which can also be found in his students and postdocs he has trained over the years: Avi Marchewka, Raz Kuperfman, Amit Singer, Boaz Nadler, Adi Taflia, Gilad Lerman, myself and many more.

When Zeev could travel, he spent one and then three months per year in Chicago during winter and summer respectively. Over the years, he had established a strong and deep collaboration with the physiologist and biophysicist Bob Eisenberg. Bob is an interesting character: telling people what he thinks even in the middle of seminars when he disagreed, but at the same time, he was deeply respecting for the role of modeling and applied mathematics. During many years, Zeev invited all of us to spend time with him at Rush, where Bob was directing, but could also talk with us about science: Bob was always pushing for considering electro-diffusion in physiology, which was clear at the ionic channel level. But, at higher scale, this theory is much harder to consider due to the complex non-cylindrical geometry. So we had to stay with diffusion for many years, until recently. Retrospectively, we have to admit that Bob was correct to push for including the electric field, because it is playing an interesting role in the conversion of current into voltage in neuronal physiology at the nano- to the micrometer level. Today, we are still working on these questions.

**The joy of working hard**

There was never a shortage of work: we found our problems to resolve in science and mostly in biology, not in any published papers. Indeed, we never found inspiration in papers, we only read the literature to write the introduction and to make sure that what we were doing was not already solved before. For many years, I usually brought the problem, formulated the mathematical aspect with Zeev and we involved whoever (students) was interested to help us resolving it. The goal was to make new discoveries in math and in the applied sciences. The mathematical aspect was only part of it.

The advantage we had is that after being trained in biology at UCSF, we could read any biological papers independently, so we could collect facts at high speed. We had very few competitors in modeling and analysis of cell biology for 10 years until 2010 where our works started to attract attention such as Narrow escape with Amit Singer, the field of modeling trafficking that we had initiated in 2002-2004, stochastic chemical reactions in 2005, Dire strait theory in 2012 and recently Extreme Narrow Escape 2017 and many others.

Indeed, in 2002, we worked on developing a theory of receptor trafficking and computing the time for escape of a receptor from a small hole was an attempt to formulate molecular trafficking at synapses, following the discoveries by R. Nicoll, R. Malenka, R. Malinow, R. Huganir understanding synaptic plasticity at a molecular level. This work was published in Jan 2004 in J.Stat Physics. We presented to the community a novel research area about modeling molecular trafficking and resolved analytically the question of exit time in dimension 2 from a small hole.  Later on, we published several papers about Narrow escape with A. Singer, a theory that had an impact in math, physics, chemistry, and simulation. We then made another discovery using Conformal Mapping that we used to resolve a similar question in more complicated domains with cusps, leading to the Dire strait time theory [6].

Another significant result was the discovery of the asymptotic expression for the full spectrum of a non-self-adjoint operator, by solving a differential equation hidden in the boundary layer of a limit cycle (the simplest obstruction in dimension 2 to a gradient field). These computations were important to revisit the distribution of time in characterizing the Up –state (depolarization) of neuronal cells.

In 2017 [15], Zeev decided to give a general talk about his road to applied mathematics at a moment of time, when applied math in Israel was suffering from emigration of young applied mathematicians, a fact that may have consequences for academic research and the high-tech economy in Israel.

**How Zeev was answering reviewers:**

Disseminating new ideas in the form of articles is an important part of the scientific activity. With no exception, all papers with Zeev contained interesting ideas and were fundamentally new. But he sometimes had hard time to publish them and also ours. Another learning experience that Zeev provided to his students is to answer reviewers and to find arguments to refute referee's comments. Constructive and good comments are rare. Most of the time, comments are a synonym of censorship. Sometimes, we had to fight hard and spent quite a bit of time in this process of finding the most incisive arguments to avoid censorship. Most criticism were usually ideological, almost never technical: comments were: "two much math", "not enough proofs", etc…Sometimes we got also very positive supportive comments, especially from J. Stat Phys, directed by Joel (Leibowitz). This was always a supportive and friendly journal, where Zeev and all of us published many original works. Finally, once a scientist has acquired a serious reputation of originality and hard worker, everybody knows will know that. The game of preventing good works to be published has usually no scientific purpose and today with modern dissemination, this is almost over due to Arxiv and BioArXiv.

**Topics covered by Schuss'scientific activity:**
-Noise activation over potential barriers
-Atomic migration in crystals
-Fast electrons in LASER plasma physics
-Lubrication in narrow bearings
-Eigenvalues of the FPO in potential fields

-Mean lifetime of the metastable state in the Josephson junction
-Loss of lock in nonlinear filtering
-Asymptotics of master equations
-Stability of large communication networks
-Noise in multi-stable systems
-Reliability of elastic structures driven by random loads
-Large deviations for Markov jump processes
-Schrodinger's equation on a lattice with weak disorder
-Ionic channels in biological membranes
-Loss of lock in RADAR tracking with jamming and manoeuvering
-Feynman integrals and absorption in quantum measurements
-Ionic simulation in a continuum
-NMR microscopy
-Molecular and cellular biology: new mathematical problems: narrow escape

**Summary of the Schuss legacy as he taught us:**

-Be prepared not to be a victim. A victim is not a hero, but remains a victim. Heroes solve great problems and people should be trained hard in that spirit.
-Applied mathematics is about making discoveries in science and technology (finding new computations of parameters, finding a framework to compute, finding new mechanisms, new chemical principles or cellular organization, etc…), rigor comes second or third: proving a theorem is no more than an exercise in polishing diamonds.
-Political power in academia: people should do their turn and after a few years, they should go back to work in their respective field or change field. Academia is about discovering, teaching and passing the knowledge.
-Students should be trained in multiple disciplines to make new discoveries and create the new science of the future.

To conclude, Zeev showed one more time that there are no standard academic paths for a life contribution to applied mathematics. This is certainly a lesson for all scholars and all universities: students should do what they like, because nobody knows from where will come the next revolution in mathematics applied to the sciences and the crisis of COVID-19 has demonstrated how modeling and simulations are keys to influence political decisions.

Zeev, we will miss you.

**Acknowledgments**: I thank Leon Glass, Amit Singer, Boaz Nadler and Bob Eisenberg for comments and corrections on this text.

* as told by Zeev to Bob Eisenberg**References (Mathematical)**

1-Schuss, Zeev. Theory and Applications of Stochastic Differential Equations (Wiley Series in Probability and Statistics - Applied Probability and Statistics Section), 1980.

2-Schuss, Zeev. *Nonlinear filtering and optimal phase tracking*. Vol. 180. Springer Science & Business Media, 2011.

3-Schuss, Zeev. *Brownian dynamics at boundaries and interfaces*. Springer-Verlag New York, 2015.

4-Schuss, Zeev. *Theory and applications of stochastic processes: an analytical approach*. Vol. 170. Springer Science & Business Media, 2009.

5-Holcman, D., & Schuss, Z. (2015). Stochastic narrow escape in molecular and cellular biology. *Analysis and Applications, Springer Verlag, NY*.

6-Holcman, D., & Schuss, Z. (2018). Second-Order Elliptic Boundary Value Problems with a Small Leading Part. In *Asymptotics of Elliptic and Parabolic PDEs* (pp. 3-9). Springer, Cham.

**References (general):**

7-Talk: Experimentalist meets theoretician held in 2007 at the Weizmann Institute, org. D. Holcman E. Korkotian. https://www.youtube.com/watch?v=lr4GTiPVOiI&t=1616s (minute 27)

8-Piano recital for conference: Experimentalist meets theoretician held in 2007 at the Weizmann Institute, org. D. Holcman E. Korkotian.

https://www.youtube.com/watch?v=taiHOlAIhDY11

9-When regression of Europe into middle ages will happen? https://www.youtube.com/watch?v=DvLGGAep6gM

10-Interview with Zeev Schuss 2005 (Israel in the Middle east., Politics of Israeli Universities, Applied mathematics)

https://www.youtube.com/watch?v=CmhshIrBXxk&t=1018s

11-Holcman, D., & Schuss, Z. , New mathematical physics needed for life sciences - Physics Today Physics Today 69, 1, 10 (2016); https://doi.org/10.1063/PT.3.3036

12-Teheran children: https://www.jewishgen.org/yizkor/Tehran/teh000.html

https://www.haaretz.com/jewish/.premium-the-tehran-children-reach-palestine-1.5308639

https://www.haaretz.com/israel-news/culture/1.4946090

http://www.isracast.com/article.aspx?id=858%20

https://rhapsodyinbooks.wordpress.com/2009/07/19/sunday-salon-review-of-%E2%80%9Cthe-children-of-teheran%E2%80%9D-a-film-by-dalia-guttman-david-tour-and-yehuda-kaveh/

13- Luxembourg Agreement : compensation of West Germany to Israel.

https://en.wikipedia.org/wiki/Reparations_Agreement_between_Israel_and_West_Germany

14- Law suit of the Teheran children, https://www.jpost.com/Jewish-World/Jewish-News/Tehran-Children-survivors-win-suit-against-state

15-Z. Schuss, A carrier Summary

https://imudotorgdotil.files.wordpress.com/2017/01/schuss2017-4.pdf

https://www.youtube.com/watch?v=3xFZxeOweCI&t=104s

16-TV Fargo series: narrow escape problem http://www.denofgeek.com/uk/tv/fargo/50276/fargo-season-3-episode-4-review-the-narrow-escape-problem

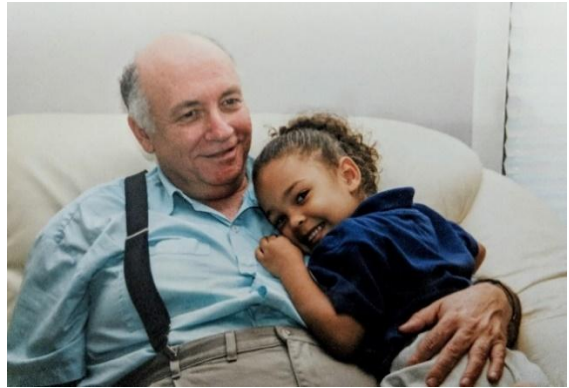

Zeev with the granddaughter of Bob Eisenberg

Poster for the 80's birthday of Z. Schuss, showing the path of the Teheran children and Zeev (left) arriving in Palestine in 1943.